\documentclass[leqno]{article}
\usepackage{amsmath}
\usepackage{amssymb}
\usepackage{amsthm}
\usepackage[usenames]{color}
\usepackage{amscd}
\usepackage{dsfont}
\usepackage[colorlinks=true,
linkcolor=webgreen,
filecolor=webbrown,
citecolor=webgreen]{hyperref}

\definecolor{webgreen}{rgb}{0,.5,0}
\definecolor{webbrown}{rgb}{.6,0,0}

\def\C{{\mathbb{C}}}

\def\N{{\mathbb{N}}}
\def\Z{{\mathbb{Z}}}

\def\1{{\bf 1}}

\def\RE{\operatorname{Re}}
\def\id{\operatorname{id}}
\def\lcm{\operatorname{lcm}}

\newtheorem{corollary}{Corollary}

\newtheorem{proposition}{Proposition}
\newtheorem{application}{Application}

\newtheorem{remark}{Remark}
\newcommand*{\wwidehat}[1]{\widehat{\smash{\widehat{#1}}\vphantom{#1}\raisebox{0.4ex}{$\vphantom{#1}$}}}

\newcommand{\DOT}{\text{\rm\Huge{.}}}

\begin{document}

\title{\bf The discrete Fourier transform of $r$-even functions}
\author{{\sc L\'aszl\'o T\'oth} \\ University of P\'ecs, Department of
Mathematics \\ Ifj\'us\'ag u. 6, H-7624 P\'ecs, Hungary \\
E-mail: ltoth@gamma.ttk.pte.hu \\  and  \\ {\sc Pentti Haukkanen} \\ Department of Mathematics 
and Statistics \\ FI-33014 University of Tampere, Finland \\ 
E-mail: pentti.haukkanen@uta.fi}
\date{}
\maketitle

\begin{abstract} We give a detailed study of the discrete Fourier transform (DFT) of $r$-even
arithmetic functions, which form a subspace of the space of
$r$-periodic arithmetic functions. We consider the DFT of sequences
of $r$-even functions, their mean values and Dirichlet series. Our
results generalize properties of the Ramanujan sum. We show that
some known properties of $r$-even functions and of the Ramanujan sum
can be obtained in a simple manner via the DFT.
\end{abstract}



\section{Introduction}

The discrete Fourier transform (DFT) of periodic functions is an important tool in various
branches of pure and applied mathematics. For instance, in number theory, the DFT of a Dirichlet
character $\chi$ (mod $r$) is the Gauss sum (character sum) given by
\begin{equation} \label{character_sum}
G(\chi,n) = \sum_{k \text{ (mod $r$)}} \chi(k) \exp(2\pi ikn/r),
\end{equation}
and if $\chi=\chi_0$ is the principal character (mod $r$), then \eqref{character_sum}
reduces to the Ramanujan sum $c_r(n)$.

For the history, properties and various applications, including signal and image processing,
of the DFT see for example the books of Briggs and Henson \cite{BriHen1995}, Broughton and Bryan
\cite{BroBry2009}, Sundararajan \cite{Sun2001}, Terras \cite{Ter1999}. For recent number
theoretical papers concerning the DFT see \cite{BecHal2010,Hau2007,Sch2008}.

It is the aim of the present paper to give a detailed study of the DFT of $r$-even arithmetic
functions, to be defined in Section 2, which form a subspace of the space of $r$-periodic arithmetic
functions.

Some aspects of the DFT of $r$-even functions were given by Haukkanen \cite{Hau2007},
Lucht \cite{Luc2000} and were considered also by Samadi, Ahmad and Swamy \cite{SamAhmSwa2005} in the
context of signal processing methods. Schramm \cite{Sch2008} investigated  the DFT of certain
special $r$-even functions, without referring to this notion.

Our results generalize and complete those of \cite{Hau2007,Luc2000,SamAhmSwa2005,Sch2008}.
Note that the Ramanujan sum $c_r(n)$ is $r$-even and it is the DFT of $\chi_0$, which is also $r$-even.
Therefore our results generalize properties of the Ramanujan sum.

The paper is organized as follows. Section \ref{section_prelim} presents an overview of the basic
notions and properties needed throughout the paper. In Section \ref{section_charact_even_func}
we give a new simple characterization of $r$-even functions. Section \ref{section_DFT_even_func}
contains properties of the DFT of $r$-even functions, while in Sections \ref{section_sequences_even_func}
and \ref{section_DFT_sequences_even_funct} we consider sequences of $r$-even functions and their
DFT, respectively. Mean values and Dirichlet series of the DFT of $r$-even functions and their sequences
are investigated in Sections \ref{section_mean_values} and \ref{section_Dirichlet_series}.

We also show that some known properties of $r$-even functions and of
the Ramanujan sum can be obtained in a simple manner via the DFT.

\section{Preliminaries} \label{section_prelim}

In this section we recall some known properties of arithmetic functions, periodic
arithmetic functions, even functions, Ramanujan sums and the DFT. We also fix the notations, most
of them being those used in the book by Schwarz and Spilker \cite{SchSpi1994}.

\subsection{Arithmetic functions} \label{subsection_arith_func}

Consider the $\C$-linear space ${\cal F}$ of arithmetic functions $f\colon\N=\{1,2,\ldots\} \to \C$
with the usual linear
operations. It is well known that with the Dirichlet convolution defined by
\begin{equation}
(f * g)(n) = \sum_{d\mid n} f(d)g(n/d)
\end{equation}
the space ${\cal F}$ forms a unital commutative $\C$-algebra. The
unity is the function $\varepsilon$ given by $\varepsilon(1)=1$ and
$\varepsilon(n)=0$ for $n>1$. The group of invertible functions is
${\cal F}^* = \{f\in {\cal F}: f(1)\ne 0\}$. The M\"obius function
$\mu$ is defined as the inverse of the function $\1 \in {\cal F}^*$
(constant $1$ function). The divisor function is $\tau= \1 * \1$,
Euler's function is $\varphi=\mu *\id$ and $\sigma=\1*\id$ is the
sum-of-divisors function, where $\id(n)=n$ ($n\in \N$). A function
$f\in {\cal F}$ is called multiplicative if $f(1)=1$ and
$f(mn)=f(m)f(n)$ for any $m,n\in \N$ such that $\gcd(m,n)=1$. The
set ${\cal M}$ of multiplicative functions is a subgroup of ${\cal
F^*}$ with respect to the Dirichlet convolution. Note that
$\1,\id,\mu,\tau,\sigma,\varphi \in {\cal M}$. For an $f\in {\cal
F}$ we will use the notation $f'=\mu*f$.

\subsection{Periodic functions} \label{subsection_periodic_func}

A function $f\in {\cal F}$ is called $r$-periodic if $f(n+r)=f(n)$
for every $n\in \N$, where $r\in \N$ is a fixed number (this periodicity extends $f$ to a function
defined on $\Z$). The set ${\cal D}_r$ of $r$-periodic functions forms an $r$-dimensional subspace of
${\cal F}$. A function $f\in {\cal F}$ is called periodic if $f\in \bigcup_{r\in \N} {\cal D}_r$.
The functions $\delta_k$ with $1\le k\le r$ given by $\delta_k(n) = 1$ for $n\equiv k$
(mod $r$) and $\delta_k(n) = 0$ for $n \not\equiv k$ (mod $r$) form a basis of ${\cal D}_r$
(standard basis).

The functions $e_k$ with $1\le k\le r$ defined by $e_k(n)=\exp(2\pi ikn/r)$ (additive characters)
form another basis of the space ${\cal D}_r$. Therefore every $r$-periodic function $f$ has
a Fourier expansion of the form
\begin{equation} \label{Fourier_exp}
f(n)= \sum_{k \text{ (mod $r$)}} g(k) \exp(2\pi i kn/r) \qquad (n\in
\N),
\end{equation}
where the Fourier coefficients $g(k)$ are uniquely determined and are given by
\begin{equation} \label{Fourier_coeff}
g(n)= \frac1{r} \sum_{k \text{ (mod $r$)}} f(k) \exp(-2\pi i kn/r)
\qquad (n\in \N)
\end{equation}
and the function $g$ is also $r$-periodic.

For a function $f\in {\cal D}_r$ its discrete (finite) Fourier transform (DFT) is the function
$\widehat{f}\in {\cal F}$ defined by
\begin{equation} \label{DFT}
\widehat{f}(n)= \sum_{k \text{ (mod $r$)}} f(k) \exp(-2\pi ikn/r) \qquad (n\in \N),
\end{equation}
where by \eqref{DFT} and \eqref{Fourier_coeff} one has $\widehat{f}=rg$.

For any $r\in \N$ the DFT is an automorphism of ${\cal D}_r$ satisfying $\skew6\wwidehat{f}=r f$.
The inverse discrete Fourier transform (IDFT) is given by
\begin{equation} \label{IDFT}
f(n)= \frac1{r}  \sum_{k \text{ (mod $r$)}} \widehat{f}(k) \exp(2\pi i kn/r)
\qquad (n\in \N).
\end{equation}

If  $f\in {\cal D}_r$, then
\begin{equation} \label{P}
\sum_{n=1}^r |\widehat{f}(n)|^2= r \sum_{n=1}^r |f(n)|^2,
\end{equation}
which is a version of Parseval's formula.

Let $f,h \in {\cal D}_r$. The Cauchy convolution 
of $f$ and $h$ is given by
\begin{equation} \label{Cauchy_periodic}
(f\otimes h)(n)= \sum_{a \text{ (mod $r$)}} f(a)h(n-a)  \qquad (n\in \N),
\end{equation}
where $({\cal D}_r,\otimes)$ is a unital commutative semigroup, the
unity being the function $\varepsilon_r$ given by
$\varepsilon_r(n)=1$ for $r\mid n$ and $\varepsilon_r(n)=0$
otherwise. Also, $\widehat{f\otimes h} = \widehat{f}\ \widehat{h}$
and $\widehat{f} \otimes \widehat{h}= r \widehat{fh}$.

For the proofs of the above statements and for further properties of $r$-periodic functions
and the DFT we refer to the books by Apostol \cite[Ch.\ 8]{Apo1976}, Montgomery and
Vaughan \cite[Ch.\ 4]{MonVau2007}, Schwarz and Spilker \cite{SchSpi1994}.

\subsection{Even functions} \label{subsection_even_func}

A function $f\in {\cal F}$ is said to be an $r$-even function if
$f(\gcd(n,r))=f(n)$ for all $n\in \N$, where $r\in \N$ is fixed.
The set ${\cal B}_r$ of $r$-even functions forms a $\tau(r)$
dimensional subspace of ${\cal D}_r$, where $\tau(r)$ is the number of
positive divisors of $r$. A function $f\in {\cal F}$ is called even if $f\in \bigcup_{r\in \N} {\cal B}_r$.
The functions $g_d$ with $d\mid r$ given by $g_d(n)=1$ if $\gcd(n,r)=d$ and $g_d(n)=0$ if
$\gcd(n,r)\ne d$ form a basis of ${\cal B}_r$. This basis can be replaced by the following one.
The functions $c_q$ with $q\mid r$ form a basis of
the subspace ${\cal B}_r$, where $c_q$ are the Ramanujan
sums, quoted in the Introduction, defined explicitly by
\begin{equation}
c_q(n)=\sum_{\substack{k \text{ (mod $q$)} \\ \gcd(k,q)=1}} \exp(2\pi ikn/q) \qquad (n,q\in \N).
\end{equation}

Consequently every $r$-even function $f$ has a (Ramanujan-)Fourier expansion of the form
\begin{equation} \label{Ramanujan_Fourier_exp}
f(n)= \sum_{d\mid r} h(d) c_d(n) \qquad (n\in \N),
\end{equation}
where the (Ramanujan-)Fourier coefficients $h(d)$ are uniquely determined and are given by
\begin{equation} \label{Ramanujan_Fourier_coeff}
h(d)= \frac1{r} \sum_{e\mid r} f(e)c_{r/e}(r/d)\qquad (d\mid r)
\end{equation}
and the function $h$ is also $r$-even. Notation: $h(d)=\alpha_f(d)$ ($d\mid r$). Note that $({\cal B}_r,\otimes)$ is a subsemigroup of $({\cal D}_r,\otimes)$ and $\alpha_{f\otimes h}(d)=r \alpha_f(d) \alpha_h(d)$ ($d\mid r$), cf. Application \ref{appl_Cauchy_convo}.

Recall the following properties of Ramanujan sums $c_r(n)$. They can be represented as
\begin{equation} \label{Ramanujan_repr}
c_r(n)=\sum_{d \mid \gcd(n,r)} d\mu(r/d) \qquad (n,r\in \N),
\end{equation}
and as
\begin{equation} \label{Holder}
c_r(n)= \frac{\mu(m) \varphi(r)}{\varphi(m)}, \quad \ m=r/\gcd(n,r),
\qquad (n,r\in \N),
\end{equation}
where \eqref{Holder} is H\"older's identity. It follows that
$c_r(n)=\varphi(r)$ for $r\mid n$ and $c_r(n)= \mu(r)$ for
$\gcd(n,r)=1$.

Let $\eta_r(n)=r$ if $r\mid n$ and $\eta_r(n)=0$ otherwise. For any fixed
$n\in \N$, $c_{\DOT}(n)=\mu*\eta_{\DOT}(n)$ and $r\mapsto c_r(n)$ is a multiplicative
function. On the other hand, $n\mapsto c_r(n)$ is multiplicative if and only if $\mu(r)=1$.

As it was already mentioned, $c_r(\DOT)$ is the DFT of the principal character (mod $r$) to be denoted
in what follows by $\varrho_r$ and given explicitly by $\varrho_r(n)=1$ if $\gcd(n,r)=1$ and
$\varrho_r(n)=0$ otherwise. Note that $\varrho_r=g_1$ with the notation of above (for $r$ fixed). Thus
\begin{equation} \label{Raman=DFT_varrho}
\widehat{\varrho}_r= c_r, \qquad \widehat{c}_r = r \varrho_r.
\end{equation}

The concept of $r$-even functions originates from Cohen
\cite{Coh1955} and was further studied by Cohen in subsequent papers
\cite{Coh1958e1,Coh1959e2,Coh1959e3}. General accounts of $r$-even
functions and of Ramanujan sums can be found in the books by
McCarthy \cite{McC1986}, Schwarz and Spilker \cite{SchSpi1994},
Sivaramakrishnan \cite{Siv1989}, Montgomery and Vaughan \cite[Ch.\
4]{MonVau2007}. See also the papers \cite{Hau2001,Spi1996,Tot2004}.


\section{Characterization of $r$-even functions} \label{section_charact_even_func}

For an $r\in \N$ let ${\cal B}'_r=\{f\in {\cal F}: f(n)=0$ for any $n\nmid r\}$.
We have

\begin{proposition} \label{prop:char_even_function}
Let $f\in {\cal F}$ and $f'=\mu * f$. Then the following assertions are equivalent:

i) $f\in {\cal B}_r$,

ii) $f(n)=\sum_{d\mid \gcd(n,r)} f'(d)$ \quad ($n\in \N$),

iii) $f'\in {\cal B}'_r$.
\end{proposition}

\begin{proof} If $f'\in {\cal B}'_r$, then for any $n\in \N$,
\begin{equation*}
f(n)=\sum_{d\mid n} f'(d)= \sum_{d\mid n, \, d\mid r}
f'(d)=\sum_{d\mid  \gcd(n,r)} f'(d)= (f'* \1)(\gcd(n,r))=
f(\gcd(n,r)).
\end{equation*}

This shows that iii) $\Rightarrow$ ii) $\Rightarrow$ i).

Now we show that i) $\Rightarrow$ iii). Assume that $f\in {\cal B}_r$ and $f'\not\in {\cal B}'_r$,
i.e., $f'(n)\ne 0$ for some $n\in \N$ with $n\nmid r$. Consider the minimal $n\in \N$ with this
property. Then all proper divisors $d$ of $n$ with $f'(d)\ne 0$ divide $r$ so that
\begin{equation*}
f(n)=\sum_{d\mid n} f'(d)= \sum_{d\mid \gcd(n,r)} f'(d)+ f'(n) = f(\gcd(n,r)) + f'(n)\ne f(\gcd(n,r)),
\end{equation*}
which gives $f\not\in {\cal B}_r$.
\end{proof}

\begin{remark} \label{rem:char_even_function} {\rm Let $f\in {\cal B}_r$. Assume that
$f(n)=\sum_{d\mid \gcd(n,r)} g(d)$ ($n\in \N$) for a function $g\in
{\cal F}$. Then $f=g\varepsilon_{\DOT}(r) * \1$ and $f=f'*\1$, by
Proposition \ref{prop:char_even_function}. Hence
$g\varepsilon_{\DOT}(r)=f'$ and obtain that $g(n)=f'(n)$ for any
$n\mid r$.}
\end{remark}

For $f=c_r$ (Ramanujan sum) we have by \eqref{Ramanujan_repr}, Proposition \ref{prop:char_even_function}
and Remark \ref{rem:char_even_function} the next identity, which can be shown also directly.

\begin{application} For any $n,r\in \N$,
\begin{equation}
\sum_{d\mid n} c_r(d)\mu(n/d)= \begin{cases} n\mu(r/n), & n\mid r, \\ 0, & n\nmid r.
\end{cases}
\end{equation}
\end{application}


\section{The DFT of $r$-even functions} \label{section_DFT_even_func}

We investigate in this section general properties of the DFT of
$r$-even functions.

\begin{proposition} \label{prop:DFT_even}
For each $r\in \N$ the DFT is an automorphism of ${\cal B}_r$. For any $f\in {\cal B}_r$,
\begin{equation} \label{evenDFTform}
\widehat{f}(n)= \sum_{d\mid r} f(d)c_{r/d}(n) \qquad (n\in \N)
\end{equation}
and the IDFT is given by
\begin{equation} \label{evenIDFTform}
f(n)= \frac1{r} \sum_{d\mid r} \widehat{f}(d)c_{r/d}(n) \qquad (n\in \N).
\end{equation}
\end{proposition}

\begin{proof}
By the definition of $r$-even functions and grouping the terms
according to the values $d=\gcd(k,r)$,
\begin{equation*}
\widehat{f}(n) = \sum_{d\mid r} f(d) \sum_{\substack{1\le j\le r/d\\ \gcd(j,r/d)=1}}
\exp(-2\pi ijn/(r/d)) = \sum_{d\mid r} f(d)c_{r/d}(n)
\end{equation*}
giving \eqref{evenDFTform} and also that $\widehat{f}\in {\cal B}_r$.
Now applying \eqref{evenDFTform} for $\widehat{f}$ (instead
of $f$) and using that $\skew6\wwidehat{f}=r f$ we have \eqref{evenIDFTform}.
\end{proof}

Proposition \ref{prop:DFT_even} is given by Lucht \cite[Th.\ 4]{Luc2000}. Formulas
\eqref{evenDFTform} and \eqref{evenIDFTform} are implicitly given by Haukkanen
\cite[Th.\ 3.2 and Eq.\ (9)]{Hau2007}, Samadi, Ahmad and Swamy \cite[Eq.\ (18)]{SamAhmSwa2005}
for $r$-even functions, and by Schramm \cite{Sch2008} for functions $n\mapsto F(\gcd(n,r))$,
where $F\in {\cal F}$ is arbitrary, without referring to the notion of even functions.

\begin{remark}  {\rm By Proposition \ref{prop:DFT_even}, for a function $f\in {\cal D}_r$ one
has $f\in {\cal B}_r$ if and only if $\widehat{f}\in {\cal B}_r$. This can be used to
show that a given function is $r$-even, cf. Application \ref{appl_Cauchy_convo}. Furthermore,
it follows that the Fourier coefficients $\alpha_f(d)$ of $f\in {\cal B}_r$ can be represented as
\begin{equation} \label{coeff_DFT}
\alpha_f(d)=\frac1{r}\widehat{f}(r/d) \qquad (d\mid r).
\end{equation}}
\end{remark}

\begin{corollary}
Let $f\in {\cal B}_r$. Then
\begin{equation} \label{sum_varphi}
\widehat{f}(n)= \sum_{d\mid r} f(d)\varphi(r/d) \qquad (r\mid n),
\end{equation}
\begin{equation} \label{sum_mu}
\widehat{f}(n)= \sum_{d\mid r} f(d)\mu(r/d) \qquad (\gcd(n,r)=1).
\end{equation}
\end{corollary}

\begin{corollary} \label{f_integer} If $f$ is a real (integer) valued $r$-even function, then
$\widehat{f}$ is also real (integer) valued.
\end{corollary}

\begin{proof} Use that $c_r(n)\in \Z$ for any $n,r\in \N$.
\end{proof}

\begin{corollary} \label{cor:evenrep}
Let $f$ be an $r$-even function. Then
\begin{equation} \label{DFT_repr} \widehat{f}(n) = \sum_{d \mid
\gcd(n,r)} d \, f'(r/d) \qquad (n\in \N),
\end{equation}
and $(\widehat{f})'(n)= n f'(r/n)$ for any $n\mid r$ and
$(\widehat{f})'(n)=0$ otherwise.
\end{corollary}

\begin{proof} Recall that $c_{\DOT}(n)=\mu*\eta_{\DOT}(n)$, see \eqref{Ramanujan_repr}. We obtain
$\widehat{f}(n)$  $=(f * c_{\DOT}(n))(r)$ $=(f * \mu *
\eta_{\DOT}(n))(r)$ $=(f'* \eta_{\DOT}(n))(r)$, and apply Remark
\ref{rem:char_even_function}.
\end{proof}

Note that by \eqref{DFT_repr} the DFT of any $f\in {\cal B}_r$ can
be written in the following forms:
\begin{equation} \label{convo_r}
\widehat{f}(n)= (f'* \eta_{\DOT}(n))(r),
\end{equation}
and
\begin{equation} \label{convo_n}
\widehat{f}= h *\1,
\end{equation}
where $h(n)=n f'(r/n)$ for $n\mid r$ and $h(n)=0$ otherwise.

\begin{proposition} \label{prop:sum_n}
Let $f$ be an $r$-even function. Then
\begin{equation}
\sum_{d\mid n} \widehat{f}(d)= \sum_{d\mid \gcd(n,r)} d\, f'(r/d)
\tau(n/d) \qquad (n\in \N).
\end{equation}
\end{proposition}

\begin{proof} Using \eqref{convo_n},
\begin{equation*}
\sum_{d\mid n} \widehat{f}(d)= (\widehat{f}*\1)(n)= (h*\1*\1)(n)= (h *\tau)(n)
=\sum_{d\mid n} h(d)\tau(n/d)
\end{equation*}
\begin{equation*}
= \sum_{d\mid \gcd(n,r)} d\, f'(r/d)\tau(n/d).
\end{equation*}
\end{proof}

In the special case $f=\varrho_r$ we reobtain (cf. \cite[Th.\ 1]{Apo1972} -- where $\sigma$ should be
replaced by $\tau$, \cite[p.\ 91]{McC1986}),
\begin{equation}
\sum_{d\mid n} c_r(d) = \sum_{d\mid \gcd(n,r)} d \mu(r/d) \tau(n/d) \qquad (n\in \N).
\end{equation}

The DFT can be used to obtain short direct proofs of certain known properties
for Ramanujan sums and special $r$-even functions. We give the
following examples.

\begin{application} {\rm By $\widehat{\varrho}_r=c_r$, cf. \eqref{Raman=DFT_varrho}, we obtain
$\skew6\wwidehat{\varrho}_r= r\varrho_r$. Therefore, by Proposition
\ref{prop:DFT_even},
\begin{equation}
\sum_{d\mid r} c_r(r/d)c_d(n)=\begin{cases} r, & \gcd(n,r)=1, \\ 0, &
\text{otherwise},
\end{cases}
\end{equation}
see \cite[p.\ 94]{McC1986}.}
\end{application}

\begin{application} {\rm Let $f(n)=(-1)^n$, which is $r$-even for any even number $r$. Its DFT is
\begin{equation}
\widehat{f}(n)= \sum_{k=1}^r (-1)^k \exp(-2\pi ikn/r)= \sum_{k=1}^r
(-\exp(-2\pi in/r))^k,
\end{equation}
which is $r$ for $n=r/2+mr$ ($m\in \Z$) and $0$ otherwise. Using
Proposition \ref{prop:DFT_even} we obtain for any even number $r$,
\begin{equation}
\sum_{d\mid r} (-1)^d c_{r/d}(n)=\begin{cases} r, & n\equiv r/2 \text{ (mod $r$)}, \\
0, & \text{otherwise}, \end{cases}
\end{equation}
cf. \cite[Th.\ IV]{NicVan1954}, \cite[p.\ 90]{McC1986}.}
\end{application}

\begin{application} \label{appl_Cauchy_convo} {\rm Let $f,h\in {\cal B}_r$. We show that their
Cauchy product $f\otimes h \in {\cal B}_r$ and the Fourier
coefficients of $f\otimes h$ are given by $\alpha_{f\otimes h}(d)= r
\alpha_f(d)\alpha_h(d)$ for any $d\mid r$, cf. Section
\ref{subsection_even_func}.

To obtain this use that $\widehat{(f\otimes h)}(n)= \widehat{f}(n)\widehat{h}(n)$ ($n\in \N$), valid
for functions $f,h\in {\cal D}_r$, cf. Section \ref{subsection_periodic_func}. Hence for any $n\in\N$,
\begin{equation*}
\widehat{(f\otimes h)}(\gcd(n,r))= \widehat{f}(\gcd(n,r)) \widehat{h}(\gcd(n,r))=
\widehat{f}(n,r) \widehat{h}(n,r)=\widehat{(f\otimes h)}(n),
\end{equation*}
showing that $\widehat{f\otimes h}$ is $r$-even. It follows that
$f\otimes h$ is also $r$-even. Furthermore, by \eqref{coeff_DFT},
for every $d\mid r$,
\begin{equation*}
\alpha_{f\otimes h}(d)= \frac1{r} (\widehat{f\otimes h})(r/d)=
\frac1{r} \widehat{f}(r/d) \widehat{h}(r/d)= r \alpha_f(d) \alpha_h(d).
\end{equation*}}
\end{application}

\begin{application} {\rm Let $N_r(n,k)$ denote the number of (incongruent) solutions (mod $r$) of the
congruence $x_1+\ldots +x_k \equiv n$ (mod $r$) with $\gcd(x_1,r)=\ldots =\gcd(x_k,r)=1$. Then it
is immediate from the definitions that
\begin{equation}
N_r(\DOT,k)= \underbrace{\varrho_r \otimes \cdots \otimes \varrho_r}_k.
\end{equation}

Therefore, $\widehat{N_r}(\DOT,k)=(\widehat{\varrho_r})^k= (c_r)^k$.
Now the IDFT formula \eqref{evenIDFTform} gives at once
\begin{equation}
N_r(n,k)= \frac1{r} \sum_{d\mid r} ((c_r(r/d))^kc_d(n) \qquad (n\in \N),
\end{equation}
formula which goes back to the work of H. Rademacher (1925) and A.
Brauer (1926) and has been recovered several times. See
\cite[Ch.\ 3]{McC1986}, \cite[p.\ 41]{SchSpi1994}, \cite{Spi1996}. }
\end{application}

\begin{application} {\rm We give a new proof of the following inversion formula of Cohen
\cite[Th.\ 3]{Coh1958e1}: If  $f$ and $g$ are $r$-even functions and if $f$ is defined by
\begin{equation}
f(n)=\sum_{d\mid r} g(d) c_d(n) \qquad (n\in \N),
\end{equation}
then
\begin{equation}
g(m)=\frac1{r} \sum_{d\mid r} f(r/d) c_d(n), \quad m=r/\gcd(n,r), \qquad (n\in \N).
\end{equation}

To show this consider the function $G(n)=g(r/\gcd(n,r))$ which is
also $r$-even. By Proposition \ref{prop:DFT_even},
\begin{equation}
\widehat{G}(n)=\sum_{d\mid r} G(r/d) c_d(n)= \sum_{d\mid r}
g(d)c_d(n) =f(n).
\end{equation}
Hence
\begin{equation}
rg(m)= rG(n)=\skew6\wwidehat{G}(n)
=\widehat{f}(n)= \sum_{d\mid r} f(r/d) c_d(n).
\end{equation}}
\end{application}

\begin{application} {\rm Anderson and Apostol \cite{AndApo1953} and Apostol \cite{Apo1972}
investigated properties of $r$-even functions $S_{g,h}$ given by
\begin{equation} \label{Apostol}
S_{g,h}(n) = \sum_{d \mid \gcd(n,r)} g(d) h(r/d) \qquad (n\in \N),
\end{equation}
where $g,h \in {\cal F}$ are arbitrary functions.

For $f=S_{g,h}$ we have according to \eqref{DFT_repr} and Remark
\ref{rem:char_even_function}, $f'(n)=g(n)h(r/n)$ ($n\mid r$) and
obtain at once
\begin{equation}
\widehat{S_{g,h}}(n)=\sum_{d\mid \gcd(n,r)} d f'(r/d)= \sum_{d\mid \gcd(n,r)} d g(r/d)h(d),
\end{equation}
which is proved in \cite[Th.\ 4]{AndApo1953} by other arguments.}
\end{application}

\begin{application} {\rm If $f$ is any $r$-even function, then
\begin{equation} \label{Parseval_even}
\sum_{n=1}^r |\widehat{f}(n)|^2 = r  \sum_{d\mid r} |f(d)|^2
\varphi(r/d).
\end{equation}

This follows by the Parseval formula \eqref{P} and grouping the
terms of the right hand side according to the values $\gcd(n,r)$.
For $f=\varrho_r$ we reobtain the familiar formula
\begin{equation}
\sum_{n=1}^r (c_r(n))^2 = r \varphi(r) \qquad (r\in \N).
\end{equation}}
\end{application}


\section{Sequences of $r$-even functions} \label{section_sequences_even_func}

In this section we consider sequences of functions $(f_r)_{r\in \N}$ such that $f_r\in {\cal B}_r$
for any $r\in \N$. Note that the sequence $(f_r)_{r\in \N}$ can be viewed also as a function of two variables:
$f\colon \N^2\to \C$, $f(n,r)=f_r(n)$.

We recall here the following concept: A function $f\colon\N^2 \to
\C$ of two variables is said to be multiplicative if $f(mn,rs)=
f(m,r)f(n,s)$ for every $m,n,r,s\in \N$ such that $\gcd(mr,ns)=1$.
For example, the Ramanujan sum $c(n,r)=c_r(n)$ is multiplicative
viewed as a function of two variables.

The next result includes a generalization of this property of the
Ramanujan sum.

\begin{proposition} \label{prop:f_multipl} Let $(f_r)_{r\in \N}$ be a sequence
of functions. Assume that

i) $f_r\in {\cal B}_r$ ($r\in \N$),

ii) $r \mapsto f_r(n)$ is multiplicative ($n\in \N$).

Then

1) the function $f \colon \N^2\to \C$, $f(n,r)=f_r(n)$ is
multiplicative as a function of two variables,

2) $f_r(m)f_r(n)=f_r(1)f_r(mn)$ holds for any $m,n,r\in \N$ with
$\gcd(m,n)=1$,

3) $n\mapsto f_r(n)$ is multiplicative if and only if $f_r(1)=1$.
\end{proposition}

\begin{proof} 1) For any $m,n,r,s\in \N$ such that $\gcd(mr,ns)=1$ we have by i) and ii),
\begin{equation*}
f_{rs}(mn) = f_r(mn) f_s(mn) = f_r(\gcd(mn,r)) f_s(\gcd(mn),s)
\end{equation*}
\begin{equation*}
=  f_r(\gcd(m,r)) f_s(\gcd(n,s))=  f_r(m) f_s(n).
\end{equation*}

2) By the definition of multiplicative functions of two variables $f
\colon \N^2\to \C$ it is immediate that $f(n,r)=\prod_p f(p^a,p^b)$
for $n=\prod_p p^a$, $r=\prod_p p^b$, and the given
quasi-multiplicative property is a direct consequence of this
equality.

3) Follows by 2).
\end{proof}

Part 1) of Proposition \ref{prop:f_multipl} is given also in
\cite{HauTot2010} and for parts 2) and 3) cf. \cite[Th.\
80]{Siv1989}.

We say that the sequence $(f_r)_{r\in \N}$ of functions is
completely even if there exists a function $F\in {\cal F}$ of a
single variable such that $f_r(n)=F(\gcd(n,r))$ for any $n,r\in \N$.
This concept originates from Cohen \cite{Coh1958e1} (for a function
of two integer variables $f(n,r)$ satisfying $f(n,r)=F(\gcd(n,r))$
for any $n,r\in \N$ he used the term completely $r$-even function,
which is ambiguous).

If the sequence $(f_r)_{r\in \N}$ is completely even, then $f_r\in
{\cal B}_r$ for any $r\in \N$, but the converse is not true. For
example, the Ramanujan sums $c_r(n)$ do not form a completely even
sequence. To see this, assume the contrary and let $p$ be any prime.
Then for $n=r=p$, $F(p)=c_p(p)=p-1$ and for $n=p, r=p^2$,
$F(p)=c_{p^2}(p)= -p$, a contradiction.

If $(f_r)_{r\in \N}$ is completely even, then
$f_r(n)=F(\gcd(n,r))=\sum_{d\mid \gcd(n,r)} F'(d)$ ($n,r\in \N$) and
by Remark \ref{rem:char_even_function} we have $f'_r(n)=F'(n)$ for
any $n\mid r$, where $F'=\mu*F$.


\section{The DFT of sequences of $r$-even functions}
\label{section_DFT_sequences_even_funct}

First we consider multiplicative properties of the DFT of sequences
of $r$-even functions

\begin{proposition} \label{prop:DFT_multipl} Let $(f_r)_{r\in \N}$ be a sequence
of functions. Assume that

i) $f_r\in {\cal B}_r$ ($r\in \N$),

ii) $r\mapsto f_r(n)$ is multiplicative ($n\in \N$).

Then

1) the function $r\mapsto \widehat{f}_r(n)$ is multiplicative ($n\in \N$),

2) the function $\widehat{f} \colon \N^2\to \C$,
$\widehat{f}(n,r)=\widehat{f}_r(n)$ is multiplicative as a function
of two variables,

3) $\widehat{f}_r(m)\widehat{f}_r(n) = f'_r(r)\widehat{f}_r(mn)$ holds for any $m,n,r\in \N$
with $\gcd(m,n)=1$,

4) $n\mapsto \widehat{f}_r(n)$ is multiplicative if and only if $f'_r(r)=1$.
\end{proposition}

\begin{proof} 1) Let $r,s\in \N$, $\gcd(r,s)=1$. Then, for any fixed $n\in \N$, by Proposition \ref{prop:DFT_even}
and using that $c_r(n)$ is multiplicative in $r$,
\begin{equation*}
\widehat{f}_{rs}(n) = \sum_{d\mid rs} f_{rs}(d) c_{rs/d}(n) =
\sum_{\substack{a\mid r\\ b\mid s}}  f_{rs}(ab) c_{(r/a)(s/b)}(n)
\end{equation*}
\begin{equation*}
= \sum_{\substack{a\mid r\\ b\mid s}} f_r(a)f_s(b) c_{r/a}(n)
c_{s/b}(n)= \sum_{a\mid r} f_r(a)  c_{n,r/a}(n) \sum_{b\mid s}
f_s(b) c_{s/b}(n)
\end{equation*}
\begin{equation*}
= \widehat{f}_r(n) \widehat{f}_s(n).
\end{equation*}

2), 3), 4) If $f_r\in {\cal B}_r$, then $\widehat{f}_r\in {\cal
B}_r$ ($r\in \N$) and by 1) we know that the function $r\mapsto
\widehat{f}_r(n)$ is multiplicative ($n\in \N$). Now apply
Proposition  \ref{prop:f_multipl} for the sequence
$(\widehat{f}_r)_{r\in \N}$ and use that $\widehat{f}_r(1)=f'_r(r)$.
\end{proof}

\begin{proposition} \label{sum_r} Let $(f_r)_{r\in \N}$ be a sequence of functions such that
$f_r\in {\cal B}_r$
($r\in \N$). Then
\begin{equation} \label{sum_even_r}
\sum_{d\mid r} \widehat{f}_d(n)= \sum_{d\mid \gcd(n,r)} d\, f_r(r/d)
\qquad (n,r\in \N),
\end{equation}
which is also $r$-even ($r\in \N$). Furthermore,
\begin{equation}
\sum_{d\mid n} \sum_{e\mid r} \widehat{f}_e(d)= \sum_{d\mid
\gcd(n,r)} d\, f_r(r/d) \tau(n/d) \qquad (n,r\in \N).
\end{equation}
\end{proposition}

\begin{proof} Similar to the proof of Proposition \ref{prop:sum_n}.
\end{proof}

In the special case $f_r=\varrho_r$ we reobtain the following known
identities for the Ramanujan sum:
\begin{eqnarray} \label{sum_c_d(n)}
\sum_{d\mid r} c_d(n) &=& \begin{cases}
         r, &  r\mid n, \\
         0, &  r\nmid n, \\
     \end{cases}\\
\sum_{d\mid n} \sum_{e\mid r} c_e(d) &=& \begin{cases}
\label{sum_sum_c_e(d)}
         r\, \tau(n/r), & r\mid n, \\
         0, &  r\nmid n,
     \end{cases}
\end{eqnarray}
\eqref{sum_c_d(n)} being a familiar one and for
\eqref{sum_sum_c_e(d)} see \cite[p.\ 91]{McC1986}.

Consider in what follows the DFT of completely even sequences,
defined in Section \ref{section_sequences_even_func}. Note that
formulae \eqref{evenDFTform} and \eqref{evenIDFTform} for the DFT
and IDFT, respectively of such sequences (that is, functions with
values $F(\gcd(n,r))$) were given by Schramm \cite{Sch2008}. He
considered also special cases of $F$.

\begin{corollary} \label{cor:F_multipl} Let $(f_r)_{r\in \N}$ be a sequence
of functions. Assume that

i) $(f_r)_{r\in \N}$ is completely even with $f_r(n)=F(\gcd(n,r))$
($n,r\in \N$),

ii) $F$ is multiplicative.

Then

1) the function $f \colon \N^2\to \C$, $f(n,r)=f_r(n)$ is
multiplicative in both variables, with the other variable fixed, and
is multiplicative as a function of two variables,

2) the function $r\mapsto \widehat{f}_r(n)$ is multiplicative ($n\in
\N$),

3) the function $\widehat{f} \colon \N^2\to \C$,
$\widehat{f}(n,r)=\widehat{f}_r(n)$ is multiplicative as a function
of two variables.

4) $n\mapsto \widehat{f}_r(n)$ is multiplicative if and only if
$F'(r)=1$.
\end{corollary}

\begin{proof} Follows from the definitions and from Proposition \ref{prop:DFT_multipl}.
\end{proof}

The results of Section \ref{section_DFT_even_func} can be applied
for completely even sequences.

\begin{corollary} Let $(f_r)_{r\in \N}$ be a completely even sequence with $f_r(n)=F(\gcd(n,r))$
($n,r\in \N$). Then
\begin{equation} \label{DFT_comp_even}
\widehat{f}_r(n)= \sum_{d\mid \gcd(n,r)} d \,F'(r/d) \qquad (n,r\in
\N),
\end{equation}
\begin{equation} \label{further_eq}
\sum_{d\mid r} \widehat{f}_{r/d}(d) =\sum_{e^2k=r} e\, F(k) \qquad
(r\in \N).
\end{equation}
\end{corollary}

\begin{proof} Here \eqref{DFT_comp_even} follows at once by Corollary \ref{cor:evenrep}, while
\eqref{further_eq} is a simple consequence of it.
\end{proof}

In particular, for $f_r=\varrho_r$ \eqref{further_eq} gives
\begin{equation}
\sum_{d\mid r} c_{r/d}(d)=\begin{cases} \sqrt{r}, & \text{$r$ is a square}, \\
0, &  \text{otherwise}, \end{cases}
\end{equation}
see \cite[p.\ 91]{McC1986}.

It follows from \eqref{DFT_comp_even} that the DFT of a completely
even sequence of functions is a special case of the functions
$S_{g,h}$ defined by \eqref{Apostol}, investigated by Anderson and
Apostol \cite{AndApo1953}, Apostol \cite{Apo1972}.

The example of $c_r(n)$ shows that the DFT sequence of a completely
even sequence is, in general, not completely even ($c_r(n)
=\widehat{\varrho}_r(n)$, where
$\varrho_r(n)=\varepsilon(\gcd(n,r))$).

Consider now the completely even sequence $f_r(n)=\tau(\gcd(n,r))$.
Then using \eqref{DFT_comp_even},
\begin{equation}
\widehat{f}_r(n) =\sum_{d\mid \gcd(n,r)} d (\mu*\tau)(r/d)=
\sum_{d\mid \gcd(n,r)} d = \sigma(\gcd(n,r))
\end{equation}
is completely even.

Next we characterize the completely even sequences such that their DFT is also a completely
even sequence.

\begin{proposition} \label{charact_compl_even} Let $(f_r)_{r\in \N}$ be a completely even sequence of
functions with $f_r(n)= F(\gcd(n,r))$. Then the DFT sequence
$(\widehat{f}_r)_{r\in \N}$ is completely even if and only if $F=c\,
\tau$, where $c\in \C$. In this case $\widehat{f}_r(n)=c\,
\sigma(\gcd(n,r))$.
\end{proposition}

\begin{proof} Assume that there is a function $G\in {\cal F}$ such that
$\widehat{f}_r(n) =\sum_{d\mid \gcd(n,r)} d\, F'(r/d)=
G(\gcd(n,r))$. Then for any $n=r\in \N$, $G(r)=
\widehat{f}_r(r)=\sum_{d\mid r} d\, F'(r/d)= (\id* F')(r)$, hence
$G$ has to be $G=\id*F'$. Now for $n=1$ and any $r\in \N$,
$G(1)=\widehat{f}_r(1)=F'(r)$. Denoting $G(1)=c$ we obtain that $F'$
is the constant function $c$. Therefore, $F=c\, \1 *\1=c\, \tau$.

Conversely, for $F=c\, \tau$ we have $F'=\mu* c\, \tau=c\, \1$ and
$\widehat{f}_r(n)= c\, \sum_{d\mid \gcd(n,r)} d= c\,
\sigma(\gcd(n,r))$.
\end{proof}

We now give a H\"older-type identity, see \eqref{Holder}, for the
DFT of completely even sequences, which is a special case of
\cite[Th.\ 2]{AndApo1953}, adopted to our case. We recall that a function $F\in {\cal F}$ said 
to be strongly multiplicative if $F$ is multiplicative and
$F(p^a)=F(p)$ for every prime $p$ and every $a\in \N$.

\begin{proposition} \label{prop:Holder} Let $(f_r)_{r\in \N}$ be a completely even
sequence with $f_r(n)=F(\gcd(n,r))$ ($n,r\in \N$). Suppose that

i) $F$ is strongly multiplicative, 

ii)  $F(p)\ne 1-p$ for any prime $p$.

Then
\begin{equation} \label{Holder_gen}
\widehat{f}_r(n)= \frac{(F*\mu)(m) (F*\varphi)(r)}{(F*\varphi)(m)},
\quad \ m=r/\gcd(n,r), \qquad (n,r\in \N).
\end{equation}

Furthermore, for every prime power $p^a$ ($a\in \N$),
\begin{equation}
\widehat{f}_{p^a}(n)=
\begin{cases}
p^{a-1}(p+F(p)-1), & p^a\mid n, \\
p^{a-1}(F(p)-1), & p^{a-1} \mid\mid n,\\
0, & p^{a-1}\nmid n.
\end{cases}
\end{equation}
\end{proposition}

\begin{proof} Here for any  prime $p$, $(F*\mu)(p)=F(p)-1$, $(F*\mu)(p^a)=0$ for
any $a\ge 2$ and $(F*\varphi)(p^a)=p^{a-1}(F(p)+p-1)$ for any $a\ge
1$. The function $F$ is multiplicative, thus $\widehat{f}_r(n)$ is
multiplicative in $r$, cf. Corollary \ref{cor:F_multipl}. Therefore
it is sufficient to verify the given identity for $r=p^a$, a prime
power. Consider three cases: Case 1) $p^a\mid n$, where
$\gcd(n,p^a)=p^a$; Case 2) $p^a\mid\mid n$, where
$\gcd(n,p^a)=p^{a-1}$; Case 3) $p^a\mid n$, where
$\gcd(n,p^a)=p^{\delta}$ with $\delta\le a-2$. \end{proof}

Recall that a function $f\in {\cal F}$ is said to be
semi-multiplicative if $f(m)f(n)= f(\gcd(m,n)) f(\lcm[m,n])$ for any
$m,n\in \N$. For example, $r\mapsto c_r(n)$ is semi-multiplicative
for any $n\in \N$. As a generalization of this property we have:

\begin{corollary} Let $(f_r)_{r\in \N}$ be a completely even sequence with $f_r(n)=F(\gcd(n,r))$
($n,r\in \N$) satisfying conditions i) and ii) of Proposition
\ref{prop:Holder}. Then $r\mapsto \widehat{f}_r(n)$ is
semi-multiplicative for any $n\in \N$.
\end{corollary}

\begin{proof} If $g\in {\cal F}$ is multiplicative, then it is known that for any
constant $C$ and any $r\in \N$, the function $n\mapsto C\,
g(r/\gcd(n,r))$ is semi-multiplicative, cf. \cite{Rea1966}, and
apply \eqref{Holder}.
\end{proof}


\section{Mean values of the DFT of $r$-even functions} \label{section_mean_values}

The mean value of a function $f\in {\cal F}$ is $m(f)=
\lim_{x\to\infty} \frac1{x}\sum_{n\le x} f(n)$ if this limit exists.
It is known that $\sum_{n\le x} c_r(n)= {\cal O}(1)$ for any $r>1$.
It follows from \eqref{Ramanujan_Fourier_exp} that the mean value of
any $r$-even function $f$ exists and is given by
$m(f)=\alpha_f(1)=\frac1{r}\widehat{f}(r)=\frac1{r}(f*\varphi)(r)$,
using \eqref{coeff_DFT}, \eqref{sum_varphi} (see also \cite[Prop.\
1]{Tot2004}). Therefore, if $f$ is $r$-even, then the mean value of
$\widehat{f}$ exists and is given by $m(\widehat{f})=\frac1{r}
\skew6\wwidehat{f}(r)=f(r)$. This follows also by Proposition
\ref{prop:DFT_even}. More exactly, we have

\begin{proposition} Let $f\in {\cal B}_r$ (with $r\in \N$ fixed).

i) If $x\in \N$ and $r\mid x$, then
\begin{equation}
\sum_{n=1}^x \widehat{f}(n) = f(r)x.
\end{equation}

ii) For any real $x\ge 1$,
\begin{equation}
\sum_{n\le x} \widehat{f}(n) = f(r)x + T_f(x), \quad |T_f(x)|\le
\sum_{d\mid r} d |f'(r/d)|.
\end{equation}

iii) The mean value of the DFT function $\widehat{f}$ is $f(r)$.
\end{proposition}

\begin{proof} For any $x\ge 1$, by Corollary \ref{cor:evenrep},
\begin{align*}
     \sum_{n\le x} \widehat{f}(n)
     &= \sum_{\substack{n\le x\\ d\mid \gcd(n,r)}} d\, f'(r/d)
     = \sum_{d\mid r} d\, f'(r/d) [x/d]
     = \sum_{d\mid r} d\, f'(r/d) (x/d-\{x/d\}) \\[\jot]
     &= x \sum_{d\mid r} f'(r/d) - \sum_{d\mid r} d\, f'(r/d) \{x/d\}
     = x f(r)+ T_f(x),
\end{align*}
where $T_f(x)$ is identically zero for $x\in \N$, $r\mid x$.
Furthermore, $T_f(x)={\cal O}(1)$ for $x\to \infty$.
\end{proof}

Now we generalize Ramanujan's formula
\begin{equation}
\sum_{n=1}^{\infty} \frac{c_r(n)}{n}= -\Lambda(r)  \qquad (r>1),
\end{equation}
where $\Lambda$ is the von Mangoldt function.

\begin{proposition} Let $f$ be an $r$-even function ($r\in \N$).

i) Then uniformly for $x$ and $r$,
\begin{equation}
\sum_{n\le x} \frac{\widehat{f}(n)}{n} = f(r) (\log x+C)-(f
*\Lambda)(r)+ {\cal O}\left(x^{-1} V_f(x)\right), \quad V_f(x)=
\sum_{d\mid r} d\, |f'(r/d)|,
\end{equation}
where $C$ is Euler's constant.

ii) If $f(r)=0$, then
\begin{equation}
\sum_{n=1}^{\infty} \frac{\widehat{f}(n)}{n} =
-(f * \Lambda)(r).
\end{equation}
\end{proposition}

\begin{proof} i) By Corollary \ref{cor:evenrep},
\begin{align*}
     \sum_{n\le x} \frac{\widehat{f}(n)}{n}
     &= \sum_{n\le x} \frac1{n} \sum_{d\mid (n,r)} d\, f'(r/d)
     = \sum_{d\mid r} f'(r/d) \sum_{j\le x/d} \frac1{j} \\[\jot]
     &= \sum_{d\mid r} f'(r/d) \bigl( \log(x/d)+ C + {\cal O}(d/x)\bigr) \\
     &= (\log x+C)\sum_{d\mid r} f'(r/d) - \sum_{d\mid r} f'(r/d) \log d
         + {\cal O}\biggl(x^{-1} \sum_{d\mid r} d|f'(r/d)|\biggr) \\
     &= (\log x+C) f(r) - (f*\mu *\log)(r)
         + {\cal O}\biggl(x^{-1} \sum_{d\mid r} d|f'(r/d)|\biggr).
\end{align*}

ii) Part ii) follows from i) with $x\to \infty$.
\end{proof}

\begin{remark} {\rm There is no simple general formula for $\sum_{r\le x}
\widehat{f}_r(n)$, where $n\in \N$ is fixed and $(f_r)_{r\in \N}$ is
a sequence of $r$-even functions (for example, $c_r(0)=\varphi(r)$
and $c_r(1)=\mu(r)$ have different asymptotic behaviors). For
asymptotic formulae concerning special functions of type
$\sum_{k=1}^n F(\gcd(k,n))$ see the recent papers
\cite{Bor2010,Tot2010}.}
\end{remark}


\section{Dirichlet series of the DFT of sequences of $r$-even functions} \label{section_Dirichlet_series}

We consider the Dirichlet series of the DFT of sequences
$(f_r)_{r\in \N}$ such that $f_r\in {\cal B}_r$ for any $r\in \N$.
By $\widehat{f}_r(n)=(\eta_{\DOT}(n)*\mu *f_r)(r)$, cf.
\eqref{convo_r}, we have formally,
\begin{align}
     \sum_{r=1}^{\infty} \frac{\widehat{f}_r(n)}{r^s}
     &= \sum_{r=1}^{\infty} \frac{\eta_r(n)}{r^s} \sum_{r=1}^{\infty}
         \frac{(f_r *\mu)(r)}{r^s}= \frac{\sigma_{s-1}(n)}{n^{s-1}}
         \sum_{r=1}^{\infty} \frac1{r^s} \sum_{k\ell=r} \mu(k) f_r(\ell) \\
     &= \frac{\sigma_{s-1}(n)}{n^{s-1}} \sum_{k=1}^{\infty}
         \frac{\mu(k)}{k^s} \sum_{\ell=1}^{\infty} \frac{f_{k\ell}(\ell)}{\ell^s}, \notag
\end{align}
where $\sigma_k(n)=\sum_{d\mid n} d^k$. This can be written in a
simpler form by considering the DFT of completely even sequences of functions.

\begin{proposition} Let $(f_r)_{r\in \N}$ be a completely even sequence of functions with
$f_r(n)=F(\gcd(n,r))$ and let  $a_F$ denote the absolute convergence
abscissa of the Dirichlet series of $F$. Then
\begin{equation}
\sum_{r=1}^{\infty} \frac{\widehat{f}_r(n)}{r^s}
=\frac{\sigma_{s-1}(n)} {n^{s-1}\zeta(s)} \sum_{r=1}^{\infty}
\frac{F(r)}{r^s}
\end{equation}
for any $n\in \N$, absolutely convergent for $\RE s>\max\{1, a_F\}$,
\begin{equation}
\sum_{n=1}^{\infty} \frac{\widehat{f}_r(n)}{n^s} =\zeta(s)
(F*\phi_{1-s})(r)
\end{equation}
for any $r\in \N$, absolutely convergent for $\RE s>1$, where
$\phi_k(r)=\sum_{d\mid r} d^k \mu(r/d)$ is a generalized Euler
function,
\begin{equation}
\sum_{n=1}^{\infty} \sum_{r=1}^{\infty}
\frac{\widehat{f}_r(n)}{n^sr^t}
=\frac{\zeta(s)\zeta(s+t-1)}{\zeta(t)} \sum_{n=1}^{\infty}
\frac{F(n)}{n^t}
\end{equation}
absolutely convergent for $\RE s>1$, $\RE t>\max\{1, a_F\}$.
\end{proposition}

\begin{proof} Apply \eqref{convo_r} and \eqref{convo_n}.
\end{proof}

For $F=\varepsilon$ we reobtain the known formulae for the Ramanujan
sum.

\vskip2mm {\bf Acknowledgement.} The authors thank Professor Lutz G.
Lucht for very helpful suggestions on the presentation of this
paper.



\begin{thebibliography}{99}

\bibitem{AndApo1953} {\sc D.~R.~Anderson} and {\sc T.~M.~Apostol}, The evaluation of
     Ramanujan's sums and generalizations, {\it Duke Math. J.} {\bf 20}
     (1953), 211--216.

\bibitem{Apo1972} {\sc T.~M.~Apostol}, Arithmetical properties of generalized
     Ramanujan sums, {\it Pacific J. Math.} {\bf 41} (1972), 281--293.

\bibitem{Apo1976} {\sc T.~M.~Apostol}, {\it Introduction to Analytic Number Theory},
Sprin\-ger, 1976.

\bibitem{BecHal2010} {\sc M.~Beck} and {\sc M. Halloran}, Finite trigonometric character sums via
discrete Fourier analysis, {\it Int. J. Number Theory} {\bf 6} (2010), 51--67.

\bibitem{Bor2010} {\sc O.~Bordell\`{e}s}, The composition of the gcd and certain arithmetic functions,
{\it J. Integer Sequences} {\bf 13} (2010), Article 10.7.1, 22 pp.

\bibitem{BriHen1995} {\sc W.~L.~Briggs} and {\sc V.~E.~Henson}, {\it The DFT -- An Owner's Manual
for the Discrete Fourier Transform}, Society for Industrial and Applied Mathematics (SIAM), 1995.

\bibitem{BroBry2009} {\sc S.~A.~Broughton} and {\sc K.~Bryan}, {\it Discrete Fourier Analysis and
Wavelets -- Applications to Signal and Image Processing}, John Wiley \& Sons, 2009.

\bibitem{Coh1955} {\sc E.~Cohen}, A class of arithmetical functions, {\it
     Proc. Nat. Acad. Sci. U.S.A.} {\bf 41} (1955), 939--944.

\bibitem{Coh1958e1} {\sc E.~Cohen}, Representations of even functions (mod
     $r$), I. Arithmetical identities, {\it Duke Math. J.} {\bf 25}
     (1958), 401--421.

\bibitem{Coh1959e2} {\sc E.~Cohen}, Representations of even functions (mod
     $r$), II. Cauchy products, {\it Duke Math. J.} {\bf 26} (1959),
     165--182.

\bibitem{Coh1959e3} {\sc E.~Cohen}, Representations of even functions (mod
     $r$), III. Special topics, {\it Duke Math. J.} {\bf 26} (1959),
     491--500.

\bibitem{Hau2001} {\sc P.~Haukkanen}, An elementary linear algebraic approach
     to even functions (mod $r$), {\it Nieuw Arch. Wiskd.} (5) {\bf 2}
     (2001), 29--31.

\bibitem{Hau2007} {\sc P.~Haukkanen}, Discrete Ramanujan-Fourier transform of even functions (mod $r$),
{\it Indian J. Math. Math. Sci.} {\bf 3} (2007), 75--80.

\bibitem{HauTot2010} {\sc P.~Haukkanen} and {\sc L. T\'oth}, An analogue of Ramanujan's sum with respect
to regular integers (mod $r$), submitted.

\bibitem{Luc2000} {\sc L.~G.~Lucht}, Discrete Fourier transform of periodic functions (Memo, unpublished),
Clausthal University of Technology, Clausthal, 2000.

\bibitem{McC1986} {\sc P.~J.~McCarthy}, {\it Introduction to Arithmetical
     Functions}, Uni\-ver\-si\-text, Sprin\-ger, 1986.

\bibitem{MonVau2007} {\sc H.~L.~Montgomery} and {\sc R.~C.~Vaughan}, {\it Multiplicative Number Theory, I.
Classical Theory}, Cambridge Studies in Advanced Mathematics 97, Cambridge University Press, 2007.

\bibitem{NicVan1954}  {\sc C.~A.~Nicol} and {\sc H.~S.~Vandiver}, A von Sterneck arithmetical function and restricted partitions with respect to a modulus, {\it Proc. Nat. Acad. Sci. U.S.A.} {\bf 40} (1954), 825--835.

\bibitem{Rea1966} {\sc D.~Rearick}, Semi-multiplicative functions, {\it Duke
     Math. J.} {\bf 33} (1966), 49--53.

\bibitem{SamAhmSwa2005} {\sc S.~Samadi}, {\sc M.~O.~Ahmad} and {\sc M.~N.~S.~Swamy}, Ramanujan sums and discrete
Fourier transforms, {\it IEEE Signal Processing Letters} {\bf 12} (2005), 293--296.

\bibitem{Sch2008} {\sc W.~Schramm}, The Fourier transform of functions of the greatest common divisor, {\it
Integers} {\bf 8} (2008), \#A50, 7 pp.

\bibitem{SchSpi1994} {\sc W.~Schwarz} and {\sc J.~Spilker}, {\it Arithmetical Functions},
London Mathematical Society Lecture Note Series, 184, Cambridge University Press, 1994.

\bibitem{Siv1989} {\sc R.~Sivaramakrishnan}, {\it Classical Theory of
     Arithmetic Functions}, in Monographs and Textbooks in Pure and
     Applied Mathematics, Vol.\ 126, Marcel Dekker, 1989.

\bibitem{Spi1996} {\sc J.~Spilker}, Eine einheitliche Methode zur Behandlung einer
linearen Kogruenz mit Nebenbedingungen, {\it Elem. Math.} {\bf 51} (1996), 107--116.

\bibitem{Sun2001} {\sc D.~Sundararajan}, {\it The Discrete Fourier Transform --
Theory, Algorithms and Applications}, World Scientific Publishing Co., 2001.

\bibitem{Ter1999} {\sc A.~Terras}, {\it Fourier Analysis on Finite Groups and Applications},
London Mathematical Society Student Texts, 43, Cambridge University Press, 1999.

\bibitem{Tot2004} {\sc L.~T\'oth}, Remarks on generalized Ramanujan sums and
     even functions, {\it Acta Math. Acad. Paedagog. Nyh\'azi. (N.S.)},
     electronic {\bf 20} (2004), 233--238.

\bibitem{Tot2010} {\sc L.~T\'oth}, A survey of gcd-sum functions, {\it J. Integer Sequences}
{\bf 13} (2010), Article 10.8.1, 23 pp.

\end{thebibliography}
\end{document}